\documentclass[11pt]{article}
\usepackage{amscd}
\usepackage{amsfonts}
\usepackage{amsmath}
\usepackage{amssymb}
\usepackage{amsthm}
\usepackage{bbm}
\usepackage{CJK}
\usepackage{fancyhdr}
\usepackage{graphicx}
\usepackage{hyperref}
\usepackage{indentfirst}
\usepackage{latexsym}
\usepackage{mathrsfs}
\usepackage{stmaryrd}
\usepackage[T1]{fontenc}
\usepackage[utf8]{inputenc}

\newtheorem{theorem}{Theorem}[section]

\newtheorem{definition}[theorem]{Definition}
\newtheorem{proposition}[theorem]{Proposition}
\newtheorem{example}[theorem]{Example}

\newtheorem{remark}[theorem]{Remark}

\usepackage[top=1in,bottom=1in,left=1.25in,right=1.25in]{geometry}
\def\<{\langle}
\def\>{\rangle}

\def\c{\cdot}

\def\o{\otimes}

\date{}
\begin{document}
\renewcommand{\baselinestretch}{1.2}
\renewcommand{\arraystretch}{1.0}
\title{\bf Twisted relative Rota-Baxter operators on Leibniz algebras and NS-Leibniz algebras  }
\date{}
\author{{\bf Apurba Das$^{1}$,   Shuangjian Guo$^{2}$\footnote
        { Corresponding author (Shuangjian Guo),  Email: shuangjianguo@126.com}}\\
{\small 1. Department of Mathematics and Statistics, Indian Institute of Technology} \\
{\small  Kanpur 208016, Uttar Pradesh, India}\\
{\small Email: apurbadas348@gmail.com}\\
{\small 2. School of Mathematics and Statistics, Guizhou University of Finance and Economics} \\
{\small  Guiyang  550025, P. R. of China}\\
  {\small Email: shuangjianguo@126.com}}
 \maketitle
\begin{center}
\begin{minipage}{13.cm}

{\bf \begin{center} ABSTRACT \end{center}}
In this paper, we introduce twisted relative Rota-Baxter operators on a Leibniz algebra as a generalization of twisted Poisson structures. We define the cohomology of a twisted relative Rota-Baxter operator $K$ as the Loday-Pirashvili cohomology of a certain Leibniz algebra induced by $K$ with coefficients in a suitable representation. Then we consider formal deformations of twisted relative Rota-Baxter operators from cohomological points of view.  Finally, we introduce and study NS-Leibniz algebras as the underlying structure of twisted relative Rota-Baxter operators.

 \medskip

 \medskip


 \medskip

 \end{minipage}
 \end{center}
 \normalsize\vskip0.5cm

\section*{Introduction}
\def\theequation{\arabic{section}. \arabic{equation}}
\setcounter{equation} {0}

The notion of Rota-Baxter operators on associative algebras was introduced in 1960 by Baxter \cite{B60} in his study of fluctuation theory in probability. Recently, it has been found many connections with dendriform algebras, pre-Lie algebras, and have applications including in Connes-Kreimer's algebraic approach to the renormalization in perturbative quantum field theory \cite{CK00}. Relative Rota-Baxter operators on Leibniz algebras were studied in \cite{ST19} which is the main ingredient in the study of the twisting theory and the bialgebra theory for Leibniz algebras. More precisely, let $(\mathfrak{g}, [\c, \c])$ be a Leibniz algebra and $(V, \rho^L, \rho^R)$ be a representation of it. A linear map $K : V \rightarrow \mathfrak{g}$ is a relative Rota-Baxter operator on $\mathfrak{g}$ with respect to the representation $V$ if it satisfies
\begin{align}\label{relative-rota-identity}
    [Ku, Kv] = K ( \rho^L (Ku) v + \rho^R (Kv) u), ~ \text{ for } u, v \in V.
\end{align}
Such operators can be seen as the Leibniz algebraic analogue of Poisson structures. Generally, Rota-Baxter operators can be defined on algebraic operads, which give rise to the splitting of operads \cite{BBGN13, PBG17}. For further details on Rota-Baxter operators, see \cite{G12}.

\medskip

Deformation theory of algebraic structures began with the seminal work of Gerstenhaber \cite{G64} for associative algebras and followed by its extension to Lie algebras by Nijenhuis and Richardson \cite{NR66, NR68}.  In general, deformation theory was developed for algebras over binary quadratic operads by Balavoine \cite{B97}. Recently, deformations of relative Rota-Baxter operators (also called $\mathcal{O}$-operators) are developed in \cite{TBGS19, Das, TSZ20}.

\medskip

In \cite{SW01} \u{S}evera and Weinstein introduced a notion of twisted Poisson structure as a Dirac structure in a certain twisted Courant algebroid. Twisted Poisson structures also studied by Klim\u{c}\'{i}k and Strobl from geometric points of view \cite{KS02}. The corresponding algebraic notion, called twisted Rota-Baxter operators was introduced by Uchino \cite{U08} in the context of associative algebras and find relations with NS-algebras of Leroux \cite{L04}. Recently, one of the present authors introduces twisted Rota-Baxter operators on Lie algebras and considers NS-Lie algebras that are related to twisted Rota-Baxter operators in the same way pre-Lie algebras are related to Rota-Baxter operators \cite{Da20}.

\medskip

Our aim in this paper is to consider twisted (relative) Rota-Baxter operators on Leibniz algebras. We show that a twisted relative Rota-Baxter operator $K$ induces a new Leibniz algebra structure and there is a suitable representation of it. The corresponding Loday-Pirashvili cohomology is called the cohomology of the twisted relative Rota-Baxter operator. As an application of the cohomology, we study deformations of a twisted relative Rota-Baxter operator $K$. We show that the infinitesimal in a formal deformation of $K$ is a $1$-cocycle in the cohomology of $K$.  Moreover, we define a notion of equivalence between two formal deformations of $K$. The infinitesimals corresponding to equivalent deformations are shown to be cohomologous. We introduce Nijenhuis elements associated with a twisted relative Rota-Baxter operator that are obtained from trivial linear deformations. We also find a sufficient condition for the rigidity of a twisted relative Rota-Baxter operator in terms of Nijenhuis elements.

\medskip

  In the last, we introduce a new algebraic structure, called NS-Leibniz algebras. We show that Ns-Leibniz algebras split Leibniz algebras and the underlying structure of a twisted relative Rota-Baxter operator. NS-Leibniz algebras also arise from Nijenhuis operators on Leibniz algebras. Further study on NS-Leibniz algebras is postponed to a forthcoming article.

\medskip

The paper is organized as follows. In Section 2, we introduce twisted relative Rota-Baxter operators in the context of Leibniz algebras and give a characterization and some new constructions. In Section 3, we define cohomology of a twisted relative Rota-Baxter operator. This cohomology has been used in Section 4 to study deformations of a twisted relative Rota-Baxter operator. Finally, in Section 5,  we introduce NS-Leibniz algebras and find its relation with twisted relative Rota-Baxter operators.

\medskip

Throughout this paper, all vector spaces, (multi)linear maps are over the field $\mathbb{C}$ of complex numbers and all the vector spaces are finite-dimensional.

\section{Leibniz algebras and Loday-Pirashvili cohomology}
\def\theequation{\arabic{section}.\arabic{equation}}
\setcounter{equation} {0}
In this section, we recall some basic definitions about Leibniz algebras and their cohomology \cite{LP93}.

\begin{definition}
A Leibniz algebra is  a vector space $\mathfrak{g}$ together with a bilinear operation (called bracket) $[\c, \c] : \mathfrak{g} \o \mathfrak{g} \rightarrow  \mathfrak{g}$  satisfying
\begin{eqnarray*}
[x, [y, z]]=[[x, y], z]+[y, [x, z]],~ \text{ for } x, y, z \in \mathfrak{g}.
\end{eqnarray*}
\end{definition}

A Leibniz algebra as above may be denoted by the pair $(\mathfrak{g}, [\c, \c])$ or simply by $\mathfrak{g}$ when no confusion arises. A Leibniz algebra whose bilinear bracket is skewsymmetric is nothing but a Lie algebra. Thus, Leibniz algebras are the non-skewsymmetric analogue of Lie algebras.

\begin{definition}
A representation of a Leibniz algebra $(\mathfrak{g}, [\c, \c])$ consists of a triple $(V, \rho^L, \rho^R)$ of a vector space $V$ and two linear maps $\rho^L, \rho^R : \mathfrak{g} \rightarrow gl(V)$ satisfying for $x, y \in \mathfrak{g}$,
\begin{align*}
\begin{cases}
\rho^L ([x,y]) = \rho^L (x) \circ \rho^L (y) - \rho^L (y) \circ \rho^L (x),\\
\rho^R ([x,y]) = \rho^L (x) \circ \rho^R (y) - \rho^R (y) \circ \rho^L (x),\\
\rho^R ([x,y]) = \rho^L (x) \circ \rho^R (y) + \rho^R (y) \circ \rho^R (x).
\end{cases}
\end{align*}
\end{definition}

It follows that any Leibniz algebra $\mathfrak{g}$ is a representation of itself with
\begin{align*}
    \rho^L (x) = L_x = [x, \c] ~~ \text{ and } ~~ \rho^R (x) = R_x = [\c, x], \text{ for } x \in \mathfrak{g}.
\end{align*}
Here $L_x$ and $R_x$ denotes the left and right multiplications by $x$, respectively. This is called the regular representation.

\medskip

Let $(\mathfrak{g}, [\c, \c])$ be a Leibniz algebra and  $(V, \rho^L, \rho^R)$ be a representation of it. The Loday-Pirashvili cohomology of $\mathfrak{g}$ with coefficients in $V$ is the
cohomology of the cochain complex $\{ C^*(\mathfrak{g}, V), \partial \}$, where $C^n (\mathfrak{g}, V) = Hom( \mathfrak{g}^{\otimes n}, V ), (n \geq 0)$ and the
coboundary operator $\partial:C^n(\mathfrak{g}, V ) \rightarrow C^{n+1}(\mathfrak{g}, V )$ given by
\begin{eqnarray*}
&&(\partial f)(x_1,\c \c \c, x_{n+1})\\
&=&\sum^{n}_{i=1}(-1)^{i+1}\rho^L(x_i)f(x_1,\c \c \c, \hat{x}_i, \c \c \c, x_{n+1}) + (-1)^{n+1}\rho^R(x_{n+1})f(x_1,\c \c \c,  x_{n})\\
&&+\sum_{1\leq i< j\leq n+1} (-1)^{i}f(x_1,\c \c \c, \hat{x}_i, \c \c \c, x_{j-1}, [x_i, x_j], x_{j+1},\c \c \c, x_{n+1}),
\end{eqnarray*}
for $x_1, \c \c \c, x_{n+1}\in \mathfrak{g}$. The corresponding cohomology groups are denoted by $H^*(\mathfrak{g}, V).$

%\begin{definition}
%An $L_{\infty}$ is a graded vector space $L=\oplus_{k}L_k$ equipped with a collection $(k \geq 1)$ of linear maps $\{l_k: \otimes^kL \rightarrow L|deg(l_k)=2-k\}$  satisfy, for any homogeneous elements $x_1, \c \c \c , x_n \in L$, \\
%(i) for any  $\sigma\in \mathbb{S}_n$,
%\begin{eqnarray*}
%l_n(x_{\sigma(1)}, \c \c \c, x_{\sigma(n-1)}, x_{\sigma(n)})=\varepsilon(\sigma) l_n(x_{1}, \c \c \c, x_{n-1}, x_{n}),
%\end{eqnarray*}
%$(ii) for $n\geq 1$,
%\begin{eqnarray*}
%\sum_{i=1}\sum_{\sigma\in \mathbb{S}_{(i, n-i)}} \varepsilon(\sigma)l_{n-i+1}(l_i(x_{\sigma(1)}, \c \c \c, x_{\sigma(i)}), x_{\sigma(i+1)}, \c \c \c , x_{\sigma(n)})=0.
%\end{eqnarray*}
%\end{definition}
%An element $\alpha \in L_1$ is said to be a Maurer-Cartan element  if it satisfies
%\begin{eqnarray*}
%l_1(\a)+\frac{1}{2!} l_2(\a, \a)-\frac{1}{3!} l_3(\a, \a, \a)-\c \c \c=0.
%\end{eqnarray*}

\section{Twisted relative Rota-Baxter operators}

In this section, we introduce twisted relative Rota-Baxter operators on Leibniz algebras and provide some examples.

\def\theequation{\arabic{section}. \arabic{equation}}
\setcounter{equation} {0}

Let $(\mathfrak{g}, [\c, \c])$ be a Leibniz algebra and $(V, \rho^L, \rho^R)$ be a representation of it.
Suppose $H \in C^2 (\mathfrak{g},V)$ is a $2$-cocycle in the Loday-Pirashvili cochain complex, i.e., $H : \mathfrak{g} \otimes \mathfrak{g} \rightarrow  V$ is a  bilinear map satisfying
\begin{small}
\begin{eqnarray*}
\rho^{L}(x)H(y, z) - \rho^{L}(y)H(x, z) - \rho^{R}(z)H(x, y) - H([x, y], z) - H(y, [x, z]) + H(x, [y, z])=0,
\end{eqnarray*}
for $x, y, z \in \mathfrak{g}.$
\end{small}

\begin{definition}
A linear map $K : V \rightarrow \mathfrak{g}$ is said to a $H$-twisted relative Rota-Baxter operator if $K$ satisfies
\begin{eqnarray*}
[Ku, Kv]=K(\rho^{L}(Ku)v+\rho^{R}(Kv)u+H(Ku, Kv)),~ \text{ for } u, v\in V.
\end{eqnarray*}
\end{definition}

\begin{example}
Any relative Rota-Baxter operator (\ref{relative-rota-identity}) is a $H$-twisted relative Rota-Baxter operator with $H = 0$.
\end{example}

\begin{example}
Let $(V, \rho^L, \rho^R)$ be a representation of a Leibniz algebra
$(\mathfrak{g}, [\c, \c])$.  Suppose $h \in C^1 (\mathfrak{g}, V)$ is an invertible $1$-cochain
in the Loday-Pirashvili cochain complex of $g$ with coefficients in $V$. Take $H = - \partial h$. Then
\begin{eqnarray*}
H(Ku, Kv) = (- \partial h)(Ku, Kv) = - \rho^{L}(K(u))v - \rho^{R}(K(v))u + h([Ku, Kv]),  ~\text{ for } u, v\in V.
\end{eqnarray*}
This shows that $K = h^{-1}: V \rightarrow  \mathfrak{g}$ is a $H$-twisted relative Rota-Baxter operator.
\end{example}

%\begin{example} {\underline{Not Correct may be}}
%Let $(\mathfrak{g}, [\c, \c])$ be  a Leibniz algebra. We denote the Leibniz bracket by the map $\mu : \mathfrak{g}^{\otimes 2} \rightarrow \mathfrak{g}$. Note that the space $V =  \mathfrak{g}^{\otimes 2}$ carries a representation of the Leibniz algebra $\mathfrak{g}$ given by
%\begin{align*}
%\rho^{L}(x) (y \otimes z) &= [x, y] \otimes z + y \otimes [x, z],\\
%\rho^{R}(x) (y \otimes z) &= y \otimes [z, x] - [x, y] \otimes z,
%\end{align*}
%for $x \in \mathfrak{g},~ y \otimes z \in  V = \mathfrak{g}^{\otimes 2}$. It is also easy to see that the map $H: \mathfrak{g}^{\otimes 2} \rightarrow \mathfrak{g}^{\otimes 2},~ y \otimes z \mapsto - y \otimes z$ is a $2$-cocycle in the Loday-Pirashvili cochain complex of $\mathfrak{g}$ with coefficients in $\mathfrak{g}^{\otimes 2}$. Then $\mu : \mathfrak{g}^{\otimes 2} \rightarrow \mathfrak{g}$ is a $H$-twisted relative Rota-Baxter operator.
%\end{example}

\begin{example}
Let $(\mathfrak{g}, [\c, \c])$ be a Leibniz algebra and $N: \mathfrak{g} \rightarrow \mathfrak{g}$ be a Nijenhuis operator it, i.e., $N$ satisfies
\begin{eqnarray*}
[Nx, Ny]=N([Nx, y] + [x,Ny] - N[x, y]), ~\text{ for } x, y \in \mathfrak{g}.
\end{eqnarray*}
In this case the vector space $\mathfrak{g}$ carries a new Leibniz algebra structure with deformed bracket
\begin{align}\label{deformed-leib}
[x, y]_N = [Nx, y] + [x,Ny] - N[x, y], \text{ for } x, y \in \mathfrak{g}.
\end{align}
This deformed Leibniz algebra $\mathfrak{g}_N = (\mathfrak{g}, [\c,\c]_N)$ has a representation on $\mathfrak{g}$ by $\rho^L(x)y:=[Nx, y]$ and $\rho^R(x)y:=[y, Nx]$, for $x \in \mathfrak{g}_N, y \in \mathfrak{g}$. With this representation, the map $H: (\mathfrak{g}_N)^{\otimes 2} \rightarrow \mathfrak{g},~H(x, y) = -N[x, y]$ is a $2$-cocycle in the Loday-Pirashvili cohomology of $\mathfrak{g}_N$ with coefficients in $\mathfrak{g}$. Moreover the identity map $Id: \mathfrak{g} \rightarrow \mathfrak{g}_N$ is a $H$-twisted relative Rota-Baxter operator.
\end{example}

\medskip

Let $K : V\rightarrow  \mathfrak{g}$ be a $H$-twisted relative Rota-Baxter operator.
 Suppose $(V', \rho'^L, \rho'^R)$ is a representation of another Leibniz algebra $(\mathfrak{g}', [\c, \c]')$  and $H' \in C^2(\mathfrak{g}', V')$ is a $2$-cocycle. Let $K' : V'\rightarrow  \mathfrak{g}'$ be a $H'$-twisted  relative Rota-Baxter operator.

\begin{definition}
A morphism of twisted relative Rota-Baxter operators from $K$ to $K'$ consists of a pair $(\phi, \psi)$ of a Leibniz algebra morphism $\phi: \mathfrak{g} \rightarrow \mathfrak{g}'$ and a linear map $\psi: V \rightarrow V'$ satisfying
\begin{align*}
&\phi\circ K=K'\circ \psi,\\
&\psi(\rho^{L}(x)u)=\rho'^{L}(\phi(x))\psi(u),~~~~\psi(\rho^{R}(x)u)=\rho'^{R}(\phi(x))\psi(u),\\
&\psi\circ H = H'\circ(\phi\o \phi), \text{ for } x\in \mathfrak{g}, u \in V.
\end{align*}
\end{definition}

\medskip

Given a $2$-cocycle $H$ in the Loday-Pirashvili cochain complex of $\mathfrak{g}$ with coefficients in $V$, one can construct the twisted semidirect product algebra. More precisely, the direct sum $\mathfrak{g} \oplus V$ carries a Leibniz algebra structure with the bilinear bracket given by
\begin{eqnarray*}
[(x,u), (y,v)]_{H} = ([x, y] ,\rho^{L}(x) v + \rho^{R}(y) u + H(x, y) ), ~ \text{ for } x, y\in \mathfrak{g}, u, v \in V.
\end{eqnarray*}
We denote this $H$-twisted semidirect product Leibniz algebra by $\mathfrak{g} \ltimes_H V$. Using this twisted semidirect product, one can characterize twisted relative Rota-Baxter operators by their graph.

\begin{proposition}\label{graph-twisted}
A linear map $K : V\rightarrow  \mathfrak{g}$ is a $H$-twisted relative Rota-Baxter operator if and only if its graph
$Gr(K ) = \{ (Ku,u)|~ u \in V\}$ is a subalgebra of the $H$-twisted semidirect product $\mathfrak{g} \ltimes_H V$.
\end{proposition}

Since $Gr (K)$ is isomorphic to $V$ as a vector space, as a consequence, we get the following.

\begin{proposition}\label{induced-leib}
Let $K : V\rightarrow  \mathfrak{g}$ be an $H$-twisted relative Rota-Baxter operator. Then the vector space $V$ carries a Leibniz algebra structure with the bracket
\begin{eqnarray*}
[u, v]_K := \rho^{L}(Ku) v + \rho^{R}(K v) u + H(K u, Kv),~ \text{ for } u, v \in  V.
\end{eqnarray*}
\end{proposition}

\subsection{Some new constructions}

In this subsection, we construct new twisted relative Rota-Baxter operators out of an old one by suitable modifications. We start with the following.

\begin{proposition}\label{prop-new-con}
Let $(V, \rho^L, \rho^R)$ be a representation of a Leibniz algebra $(\mathfrak{g}, [\c, \c])$. For any $2$-cocycle $H \in C^2(\mathfrak{g}, V)$ and $1$-cochain $h \in C^1 (\mathfrak{g}, V)$, the twisted Leibniz algebras $\mathfrak{g} \ltimes_H V$ and $\mathfrak{g} \ltimes_{H + \partial h} V$ are isomorphic.
\end{proposition}

{\bf Proof.} We define an isomorphism $\Psi_h: \mathfrak{g} \ltimes_H V \rightarrow \mathfrak{g} \ltimes_{H + \partial h} V$ of the underlying vector spaces by $\Psi_h (x,u):= (x,u - h(x))$, for $(x,u) \in \mathfrak{g} \ltimes_H V$. Moreover, we have
\begin{eqnarray*}
&& \Psi_h([(x,u), (y,v)]_{H}) \\
&=& \Psi_h([x, y],\rho^{L}(x)v+\rho^{R}(y)u+H(x, y))\\
&=& ([x, y],\rho^{L}(x)v + \rho^{R}(y)u + H(x, y) - h[x, y])\\
&=&([x, y],\rho^{L}(x)v + \rho^{R}(y)u + H(x, y) - \rho^{L}(x)h(y) - \rho^{R}(y)h(x) + \partial h(x, y) )\\
&=&[(x,u-h(x)), (y,v-h(y))]_{H +\partial h}\\
&=&[\Psi_h(x,u), \Psi_h(y,v)]_{H +\partial h}.
\end{eqnarray*}
This shows that $\Psi_h$ is in fact an isomorphism of Leibniz algebras. \hfill $\square$

\begin{proposition}
Let $K : V \rightarrow \mathfrak{g}$ be a $H$-twisted Rota-Baxter operator. For any $1$-cochain $h \in C^1(\mathfrak{g}, V)$, if the linear map $(Id_V - h \circ K) : V \rightarrow V$ is invertible then the map $K \circ (Id_V - h \circ K)^{-1} : V \rightarrow \mathfrak{g}$ is a $(H + \partial h)$-twisted relative Rota-Baxter operator.
\end{proposition}

{\bf Proof.}
Consider the subalgebra $Gr(K) \subset \mathfrak{g} \ltimes_H V$ of the $H$-twisted semidirect product. Thus by Proposition \ref{prop-new-con}, we get that
\begin{align*}
\Psi_h (Gr(K))=\{ (K u,u-hK(u))|~u \in V\} \subset \mathfrak{g} \ltimes_{H + \partial h} V
\end{align*}
is a subalgebra. Since the map $(Id_V -h \circ K ): V\rightarrow V$ is invertible, we have $\Psi_h (Gr(K))$ is the graph of the linear map $K \circ (Id_V - h \circ K)^{-1}$. Hence by Proposition \ref{graph-twisted}, the map $K \circ (Id_V - h \circ K)^{-1}$ is a $(H + \partial h)$-twisted relative Rota-Baxter operator. \hfill $\square$

\medskip

Let $K : V \rightarrow \mathfrak{g}$ be a $H$-twisted Rota-Baxter operator. Suppose $B \in C^1(\mathfrak{g}, V)$ is a $1$-cocycle in the Loday-Pirashvili cochain complex of $\mathfrak{g}$ with coefficients in $V$. Then $B$ is said to be $K$-admissible if the linear map $(Id_V + B \circ K ) : V \rightarrow V$ is invertible. With this notation, we have the following.

\begin{proposition}
Let $B \in C^1 (\mathfrak{g}, V)$ be a $K$-admissible $1$-cocycle. Then the map $K \circ (Id_V + B \circ K )^{-1}: V\rightarrow g$ is a $H$-twisted Rota-Baxter operator.
\end{proposition}

{\bf Proof.} Consider the deformed subspace
\begin{eqnarray*}
\tau_B (Gr(K))=\{ (K u, u + B\circ K(u))|~u \in V\} \subset \mathfrak{g} \ltimes_{H} V.
\end{eqnarray*}
Since $B$ is a $1$-cocycle, $\tau_B (Gr (K)) \subset \mathfrak{g} \ltimes_{H} V$  turns out to be a subalgebra. Further, the map $(Id_V + B \circ K )$ is invertible implies that $\tau_B (Gr (K))$ is the graph of the map $K \circ (Id_V + B \circ K )^{-1}$. Hence the result follows from Proposition \ref{graph-twisted}. \hfill $\square$

\medskip

The $H$-twisted Rota-Baxter operator in the above proposition is called the gauge
transformation of $K$ associated with $B$. We denote this $H$-twisted relative Rota-Baxter operator simply by $K_B$.

\begin{proposition}
Let $K$ be a $H$-twisted relative Rota-Baxter operator and $B$ be a $K$-admissible $1$-cocycle. Then
the Leibniz algebra structures on $V$ induced from the $H$-twisted Rota-Baxter operators $K$ and $K_B$ are isomorphic.
\end{proposition}

{\bf Proof.} Consider the linear isomorphism $(Id_V + B \circ K) : V \rightarrow V$. Moreover, for any $u, v \in V$, we have
\begin{eqnarray*}
&& [(Id_V + B \circ K )(u), (Id_V + B \circ K )(v)]_{K_B}\\
&=& \rho^{L}(K(u)) (Id_V + B \circ K )(v) +\rho^{R}(K (v))(Id_V + B \circ K )(u) + H(K u, Kv)\\
&=&  \rho^{L}(K(u))v+ \rho^{R}(K(v))u+ \rho^{L}(K(u)) (B \circ K(v)) + \rho^{R}(K(v)) (B \circ K(u)) + H(K u, Kv)\\
%&=& \rho^{L}(K(u))v+ \rho^{R}(K(v))u+ \rho^{L}(K(u)) B \circ K(v)+\rho^{R}(K(v)) B \circ K(u)+ H(K u, Kv)\\
&=& \rho^{L}(K(u))v + \rho^{R}(K(v))u + B[Ku, Kv] + H(K u, Kv)\\
&=& [u, v]_K+B\circ K([u, v]_K)\\
&=& (Id_V + B \circ K)([u, v]_K).
\end{eqnarray*}
This shows that $(Id_V + B \circ K) : (V, [\c, \c]_K) \rightarrow (V, [\c, \c]_{K_B})$ is a Leibniz algebra isomorphism.  \hfill $\square$

\section{Cohomology of twisted relative Rota-Baxter operators}

In this section, we define cohomology of a $H$-twisted relative Rota-Baxter operator $K$ as the Loday-Pirashvili cohomology of the Leibniz algebra $(V, [ \c,\c ]_K )$ constructed in Proposition \ref{induced-leib} with coefficients in a suitable representation on $\mathfrak{g}$. In the next section, we will use this cohomology to study deformations of $K$.

\begin{proposition}
Let $K : V \rightarrow \mathfrak{g}$ be a $H$-twisted relative Rota-Baxter operator. Define maps $\overline{\rho}^L, \overline{\rho}^R : V \rightarrow gl (\mathfrak{g})$ by
\begin{small}
\begin{align*}
\overline{\rho}^L(u)x = [Ku, x] - K (\rho^R(x)u) - KH(Ku, x),~~~ \overline{\rho}^R(u)x = [x, Ku] - K(\rho^L(x)u)-KH(x, Ku),
\end{align*}
\end{small}
for $ u\in V$ and $ x\in \mathfrak{g}$. Then $(\mathfrak{g},  \overline{\rho}^L, \overline{\rho}^R)$ is a representation of the Leibniz algebra $(V, [\c, \c]_{K})$.
\end{proposition}

{\bf Proof.}  For $u, v \in V$ and $x \in \mathfrak{g}$ , we have
\begin{small}
\begin{eqnarray*}
&& \overline{\rho}^L(u)\overline{\rho}^L(v)x-\overline{\rho}^L(v)\overline{\rho}^L(u)x\\
&=& \overline{\rho}^L(u)([Kv, x] - K(\rho^R(x)v) - KH(Kv, x)) - \overline{\rho}^L(v)([Ku, x] - K(\rho^R(x)u) - KH(Ku, x))\\
%&=& [Ku, [Kv, x] - K(\rho^R(x)v) - KH(Kv, x)] - K (\rho^R([Kv, x])u) + K(\rho^R(K\rho^R(x)v)u)\\
%&&+K (\rho^R(KH(Kv, x))u) - KH(Ku, [Kv, x] - K(\rho^R(x)v)-KH(Kv, x))\\
%&&-[Kv, [Ku, x] - K(\rho^R(x)u) - KH( Ku, x)] + K(\rho^R([Ku, x])v) - K (\rho^R(K\rho^R(x)u)v)\\
%&&-K (\rho^R(KH(Ku, x))v)+KH(Kv, [Ku, x]-K(\rho^R(x)u)-KH(x, Ku))\\
&=&[Ku, [Kv, x]] - [Ku,K(\rho^R(x)v)] - [Ku,KH(Kv, x)] - K (\rho^R([Kv, x])u)
+K (\rho^R(K\rho^R(x)v)u) \\
&&+ K (\rho^R(KH( Kv, x))u) - KH(Ku, [Kv, x])  + KH(Ku,  K(\rho^R(x)v)) + KH(Ku, KH( Kv, x))\\
&&-[Kv, [Ku, x]] + [Kv, K(\rho^R(x)u)] + [Kv, KH(Ku, x)] + K(\rho^R([Ku, x])v)
-K(\rho^R(K\rho^R(x)u)v) \\ &&- K(\rho^R(KH(Ku, x))v) + KH(Kv, [Ku, x]) - KH(Kv, K(\rho^R(x)u)) - KH(Kv, KH(Ku, x))\\
%&=&[Ku, [Kv, x]] - [Ku,K(\rho^R(x)v)] - K (\rho^L(Kv)\rho^R(x)u)\\
%&&+K(\rho^R(x)\rho^L(Kv)u) + K(\rho^R(K\rho^R(x)v)u) - [Kv, [Ku, x]]\\
%&&+[Kv, K(\rho^R(x)u)] + K(\rho^L(Ku)\rho^R(x)v) - K(\rho^R(x)\rho^L(Ku)v)\\
%&& -K(\rho^R(K\rho^R(x)u)v)-KH(K[u, v]_{K}, x)\\
&=& [[Ku, Kv], x] - K(\rho^R(x)\rho^L(K u)v) - K (\rho^R(x)\rho^R(K v)u) - K (\rho^R (x) H (Ku,Kv))- KH(K[u, v]_{K}, x)\\
&=&[K[u, v]_{K}, x] - K (\rho^R(x)[u, v]_{K}) - KH(K[u, v]_{K}, x)\\
&=&\overline{\rho}^L([u, v]_{K})x.
\end{eqnarray*}
\end{small}
\noindent The third equality is obtained by some cancellations and using the fact that $H$ is a $2$-cocycle.
Thus, we deduce that
 \begin{eqnarray*}
\overline{\rho}^L([u, v]_{K})=\overline{\rho}^L(u)\overline{\rho}^L(v) - \overline{\rho}^L(v)\overline{\rho}^L(u).
 \end{eqnarray*}
We also have
\begin{small}
\begin{eqnarray*}
&&  \overline{\rho}^L(u)\overline{\rho}^R(v)x-\overline{\rho}^R(v)\overline{\rho}^L(u)x\\
&=& \overline{\rho}^L(u)([x, Kv] - K(\rho^L(x)v)-KH(x, Kv))-\overline{\rho}^R(v)([Ku, x] - K(\rho^R(x)u) - KH( Ku, x))\\
%&=& [Ku, [ x, Kv] - K (\rho^L(x)v)-KH(x, Kv)] - K\rho^R([x, Kv] -K(\rho^L(x)v)-KH(x, Kv))u\\
%&& -KH(Ku, [ x, Kv] -K(\rho^L(x)v)-KH(x, Kv))-[[Ku, x]- K(\rho^R(x)u)-KH( Ku, x), Kv]\\
%&&+K \rho^L([Ku, x]- K(\rho^R(x)u) - KH( Ku, x))v + KH([Ku, x] - K(\rho^R(x)u) - KH( Ku, x), Kv)
&=&[Ku, [x, Kv]] - [Ku,K(\rho^L(x)v)] - [Ku, KH(x, Kv)] - K(\rho^R([ x, Kv])u) +K (\rho^R(K\rho^L(x)v)u) \\&&+ K(\rho^R(KH(x, Kv))u) - KH(Ku, [ x, Kv]) + KH(Ku,K(\rho^L(x)v)) + KH(Ku,KH(x, Kv))\\
&& - [[Ku, x], Kv] + [ K(\rho^R(x)u), Kv] + [KH( Ku, x), Kv] + K (\rho^L([Ku, x])v) - K(\rho^L(K\rho^R(x)u)v) \\ &&- K (\rho^L(KH( Ku, x))v) + KH([Ku, x], Kv) - KH(K(\rho^R(x)u), Kv) - KH( KH( Ku, x), Kv)\\
%&=& [Ku, [x, Kv]] - [Ku,K(\rho^L(x)v)] - K(\rho^L(x)\rho^R(Kv)u)\\
%&& + K(\rho^R(Kv)\rho^L(x)u) + K(\rho^R(K\rho^L(x)v)u) - [[Ku, x], Kv]\\
%&& + [ K(\rho^R(x)u), Kv] + K(\rho^L(Ku)\rho^L(x)v) - K (\rho^L(x)\rho^L(Ku)v)\\
%&& - K (\rho^L(K\rho^R(x)u)v) - KH(x, K[u, v]_{K})\\
&=&[x, [Ku,Kv]] - K (\rho^L (x) \rho^L (Ku)v) -  K (\rho^L (x) \rho^R (Kv) u) - K\rho^L (x) H(Ku, Kv))  - KH(x, K[u, v]_{K})\\
&=&[x, K[u, v]_K] - K(\rho^L(x)[u,v]_K) - KH(x, K[u, v]_{K})\\
&=& \overline{\rho}^R([u, v]_{K})x
\end{eqnarray*}
\end{small}
which shows that
\begin{eqnarray*}
\overline{\rho}^R([u, v]_{K})=\overline{\rho}^L(u)\overline{\rho}^R(v)-\overline{\rho}^R(v)\overline{\rho}^L(u).
 \end{eqnarray*}
 Similarly, we can show that
 \begin{eqnarray*}
\overline{\rho}^R([u, v]_K) = \overline{\rho}^R(v) \circ \overline{\rho}^R(u) +  \overline{\rho}^L(u) \circ \overline{\rho}^R(v).
 \end{eqnarray*}
 Therefore, $(\mathfrak{g},  \overline{\rho}^L, \overline{\rho}^R)$ is a representation of the Leibniz algebra $(V, [\c, \c]_{K})$.  \hfill $\square$

 \medskip

 We will now consider the Loday-Pirashvili cohomology of the
 Leibniz algebra $(V, [\c, \c]_K)$ with coefficients in the representation $(\mathfrak{g}, \overline{\rho}^L, \overline{\rho}^R)$. More precisely, we define
\begin{eqnarray*}
C^n(V, \mathfrak{g}):= Hom(V^{\otimes n}, \mathfrak{g}),~ \text{ for } n\geq 0.
\end{eqnarray*}
and the differential $\partial_K: C^n(V, \mathfrak{g}) \rightarrow  C^{n+1}(V, \mathfrak{g})$ by
\begin{align*}
&(\partial_Kf)(u_1, \c \c \c, u_{n+1})\\
=& \sum_{i=1}^{n}(-1)^{i+1}[Ku_i, f(u_1, \c \c \c, \hat{u}_i, \c \c \c, u_{n+1})] - \sum_{i=1}^{n} (-1)^{i+1} K (\rho^R(f(u_1, \c \c \c, \hat{u}_i, \c \c \c, u_{n+1})) u_i)\\
&-\sum_{i=1}^{n}(-1)^{i+1}KH( Ku_i, f(u_1, \c \c \c, \hat{u}_i, \c \c \c, u_{n+1}))+(-1)^{n+1}[f(u_1, \c \c \c, u_{n}), Ku_{n+1}]\\
& +(-1)^{n} K (\rho^{L}(f(u_1, \c \c \c, u_{n}))u_{n+1}) +(-1)^{n}KH( f(u_1, \c \c \c, u_{n}), Ku_{n+1})\\
&+\sum_{1\leq i< j\leq n+1}(-1)^if(u_1, \c \c \c, \hat{u}_i,\c \c \c, u_{j-1},\rho^L(Ku_i)u_j +\rho^R(Ku_j)u_i+H(Ku_i, Ku_j), u_{j+1}, \c \c \c, u_{n+1}),
\end{align*}
for $f\in C^n(V, \mathfrak{g})$ and $u_1, \c \c \c, u_{n+1} \in V$. Denote by
\begin{align*}
    Z^n (V, \mathfrak{g}) =~& \{ f \in C^n (V, \mathfrak{g}) ~|~ \partial_K f = 0 \}, \\
    B^n (V, \mathfrak{g}) =~& \{ \partial_K g ~|~ g \in C^{n-1}(V, \mathfrak{g}) \}.
\end{align*}
The corresponding cohomology groups
\begin{align*}
    H^n (V, \mathfrak{g}) := \frac{  Z^n (V, \mathfrak{g})  }{B^n (V, \mathfrak{g})}, \text{ for } n \geq 0
\end{align*}
are called the cohomology of the $H$-twisted Rota-Baxter operator $K$.

%\begin{proposition}
%Let $K$ be a $H$-twisted relative Rota-Baxter operator. Then for any $f\in Hom(V^{\otimes n}, \mathfrak{g})$, we have
%\begin{eqnarray*}
%d_K f=(-1)^{n-1}\partial_K f.
%\end{eqnarray*}
%Consequently, the cohomology of $K$ is isomorphic to the Chevalley-Eilenberg cohomology of $(V, [\c , \c]_K)$ with
%coefficients in $g$.
%\end{proposition}

%{\bf Proof.} For all $u_1, u_2, \c \c \c, u_{n+1} \in  V$, one have
%\begin{eqnarray*}
%&& \llbracket K,f \rrbracket (u_1, \c \c \c, u_{n+1})\\
%&=& -K(\rho^L(f(u_1, \c \c \c, u_{n}))u_{n+1})+\sum_{i=1}^n(-1)^{n+1-i}K(\rho^R(f(u_1, \c \c \c,\hat{u}_i,\c \c \c, u_{n+1}))u_i)\\
%&&+[f(u_1, \c \c \c, u_{n}), Ku_{n+1}] +\sum_{i=1}^n(-1)^{n+i}[Ku_i, f(u_1, \c \c \c,\hat{u}_i,\c \c \c, u_{n+1})]\\
%&&+\sum_{1\leq i < j\leq n+1}(-1)^{n-1-i}f(u_1, \c \c \c,\hat{u}_i,\c \c \c, u_{j-1}, \rho^L(Ku_i)u_j + \rho^R(Ku_j)u_i, u_{j+1}, \c \c \c, u_{n+1}).
%\end{eqnarray*}
%See \cite{TSZ20} for details. Next from the expression (\ref{ternary-brkt}), we see that
%\begin{eqnarray*}
%&& \llbracket K, K, f \rrbracket (u_1, \c \c \c, u_{n+1})\\
%&=&  (-1)^{n}~ 2 \big(  \sum_{i=1}^{n}(-1)^{i} KH( Ku_i, f(u_1, \c \c \c, \hat{u}_i, \c \c \c, u_{n+1})) + (-1)^n KH (f (u_1, \c \c \c, u_n), Ku_{n+1} )\\
%&& \qquad \quad +\sum_{1\leq i< j\leq n+1}(-1)^{i}f(u_1, \c \c \c, \hat{u}_i,\c \c \c, u_{j-1}, H(Ku_i, Ku_j), u_{j+1}, \c \c \c, u_{n+1}) \big).
%\end{eqnarray*}
%Hence, $d_K f=\llbracket K,f \rrbracket-\frac{1}{2}\llbracket K, K, f \rrbracket =(-1)^{n-1}\partial_K f$. Therefore, the proof is finished. \hfill $\square$

\section{Deformations of twisted relative Rota-Baxter operators}

\def\theequation{\arabic{section}. \arabic{equation}}
\setcounter{equation} {0}

In this section, we will apply the classical deformation theory of Gerstenhaber to twisted relative Rota-Baxter operators. We will introduce certain elements (called
Nijenhuis elements) associated with a $H$-twisted relative Rota-Baxter operator that arise from trivial linear deformations. We also consider rigidity of a $H$-twisted relative Rota-Baxter operator and give a sufficient condition for rigidity in terms of Nijenhuis elements.

\subsection{Linear deformations}
Let $(\mathfrak{g}, [\c,\c])$ be a Leibniz algebra, $(V,  \rho^L, \rho^R)$ be a representation of it, and $H \in C^2 (\mathfrak{g}, V)$ be a $2$-cocycle in the Loday-Pirashvili cochain complex. Suppose $K : V \rightarrow \mathfrak{g}$ is a $H$-twisted relative Rota-Baxter operator.

\begin{definition}
A linear map $K_1 : V \rightarrow \mathfrak{g}$ is said to generate a linear deformation of the $H$-twisted relative Rota-Baxter operator $K$ if for all $t \in \mathbb{C}$, the sum $K_t = K +t K_1$ is still a $H$-twisted relative Rota-Baxter operator. In this case, $K_t = K +t K_1$ is said to be a linear deformation of $K$.
\end{definition}

Suppose $K_1$ generates a linear deformation of $K$. Then we have
\begin{align*}
    [K_t u , K_t v] = K_t \big(   \rho^L ( K_t u) v + \rho^R (K_t v) u + H (K_t u, K_tv) \big), ~ \text{for } u, v \in V.
\end{align*}
This is equivalent to the following conditions
\begin{align}
[Ku, K_1v] + [K_1u, K v] =~& K_1(\rho^L(K u)v + \rho^R(K v)u+H(Ku, Kv))\nonumber\\
&+K (\rho^L(K_1u)v + \rho^R(K_1v)u+H(K_1u, Kv)+H(Ku, K_1v)),\\
[K_1u, K_1v]
 =K_1(\rho^L(K_1u)v ~&+ \rho^R(K_1v)u+H(Ku, K_1v)+H(K_1u, Kv))+KH(K_1u, K_1v),\\
K_1(H(K_1(u), K_1(v)))=~&0.
\end{align}

Note that Eq. (4.1) means that $K_1$ is a $1$-cocycle in the cohomology of $K$. Hence $K_1$ induces an element in $H^1_K (V, \mathfrak{g}).$

\begin{definition}
Two linear deformations $K_t = K +tK_1$  and $K'_t = K +tK'_1$ of $K$ are said to be equivalent if there exists an element $x \in \mathfrak{g}$ such that
\begin{align*}
(\phi_t = Id_\mathfrak{g} +  t L_x,~ \psi_t = Id_V + t(\rho^{L}(x)+H(x, K-))
\end{align*}
is a morphism of twisted Rota-Baxter operators from $K_t$ to $K'_t$.
\end{definition}

The condition that $\phi_t = Id_\mathfrak{g} +  t L_x$ is a Leibniz algebra morphism of $(\mathfrak{g}, [\c, \c])$ is equivalent to
\begin{eqnarray}\label{alg-map-imply}
[[x, y], [x, z]]=0,~ \text{ for } y, z\in \mathfrak{g}.
\end{eqnarray}
Further, The conditions $\psi_t(\rho^{L}(y)u) = \rho^{L}(\phi_t(y)) \psi_t(u)$ and $\psi_t(\rho^{R}(y)u) = \rho^{R}(\phi_t(y)) \psi_t(u)$, for $y \in \mathfrak{g}, u \in V$ are respectively equivalent to
\begin{eqnarray}\label{left-action-imply}
\left\{
\begin{aligned}
H(x, K(\rho^{L}(y)u)) = \rho^{L}(y)H(x, Ku), \\
\rho^{L}([x, y]) (\rho^{L}(x)u+H(x, Ku))= 0,
\end{aligned}
\right. \end{eqnarray}
\begin{eqnarray}\label{right-action-imply} \left\{
\begin{aligned}
H(x, K(\rho^{R}(y)u)) = \rho^{R}(y)H(x, Ku),\\
\rho^{R}([x, y]) (\rho^{L}(x)u+H(x, Ku))= 0.
\end{aligned}
\right. \end{eqnarray}
Similarly, the conditions $\psi_t \circ H = H \circ (\phi_t \o \phi_t)$ and $\phi_t \circ K_t = K'_t \circ \psi_t$ are respectively equivalent to
\begin{eqnarray}\label{h-comp-imply}
\left\{
\begin{aligned}
&\rho^{L}(x)H(y, z)+H(x, KH(y, z)) = H(x, [y, z])+H(y, [x, z]),\\
&H([x, y], [x, z])= 0,
\end{aligned}
\right. \end{eqnarray}
\begin{eqnarray}\label{diffe}
 \left\{
\begin{aligned}
& K_1(u)+[x, Ku]=K(\rho^{L}(x)u+H(x, Ku)) + K'_1(u), \\
&[x, K_1u] = K'_1 (\rho^{L}(x)u+H(x, Ku)).
\end{aligned}
\right.
\end{eqnarray}

It follows from the first identity in (\ref{diffe}) that  $K_1(u)-K'_1(u) = d_{K}(x)(u)$. Hence we obtain the following.

\begin{theorem}
If two linear deformations
$K_t = K +tK_1$  and $K'_t = K +tK'_1$ of a $H$-twisted relative Rota-Baxter operator $K$ are equivalent, then $K_1$ and $K'_1$ are in the same
cohomology class of $H^1_K (V, \mathfrak{g})$.
\end{theorem}

\begin{definition}
 A linear deformation $K_t = K +tK_1$ of a $H$-twisted relative Rota-Baxter operator $K$ is said to be trivial if
$K_t$ is equivalent to the undeformed deformation $K'_t=K$.
\end{definition}

We will now define Nijenhuis elements associated with a $H$-twisted relative Rota-Baxter operator $K$ in a way that a trivial deformation of $K$ induces a Nijenhuis element.

\begin{definition}
 Let $K$ be a $H$-twisted relative Rota-Baxter operator. An element $x\in \mathfrak{g}$ is called a Nijenhuis
element associated with $K$ if $x$ satisfies
\begin{eqnarray*}
[x, \overline{\rho}^{R}(u)(x)]=0, ~\text{ for } u \in V
\end{eqnarray*}
and Equations (\ref{alg-map-imply}), (\ref{left-action-imply}), (\ref{right-action-imply}), (\ref{h-comp-imply}) hold.
\end{definition}

The set of all Nijenhuis elements associated with $K$ is denoted by Nij$(K)$. As mentioned earlier that a trivial deformation induces a Nijenhuis element. In the next subsection, we give a sufficient condition for the rigidity of a twisted Rota-Baxter operator in terms of Nijenhuis elements.

\subsection{Formal deformations }

%Let $(\mathfrak{g}, [\c, \c])$ be a Leibniz algebra,  $(V,  \rho^L, \rho^R)$ be a representation of it, and $H \in C^2 (\mathfrak{g}, V)$ be a $2$-cocycle in the Loday-Pirashvili cochain complex of $\mathfrak{g}$ with coefficients in $V$. Suppose $K : V \rightarrow \mathfrak{g}$ is a $H$-twisted relative Rota-Baxter operator.

Let $\mathbb{C}[[t]]$ be the ring of power series in one variable $t$. For any $\mathbb{C}$-linear space $V$, we
let $V [[t]]$ denotes the vector space of formal power series in $t$ with coefficients in $V$. Moreover, if $(\mathfrak{g}, [\c, \c])$ is a Leibniz algebra over $\mathbb{C}$, then one can extend the Leibniz bracket on $\mathfrak{g} [[t]]$ by $\mathbb{C}[[t]]$-bilinearity.
%there is a $\mathbb{C}[[t]]$-Leibniz
%algebra structure on $g[[t]]$ given by
%\begin{eqnarray*}
%&&[\sum_{i=0}^{+\infty}x_it^{i},\sum_{j=0}^{+\infty}y_jt^{j} ]_g=\sum_{k=0}^{+\infty}\sum_{i+j=k}[x_i, y_j]t^{k}, ~~~\forall x_i, y_j\in g.
%\end{eqnarray*}
Furthermore, if $(V, \rho^L, \rho^{R})$ is a representation of the Leibniz algebra $(\mathfrak{g}, [\c, \c])$, then there is a representation $(V [[t]],  \rho^L, \rho^{R})$ of the Leibniz algebra $\mathfrak{g}[[t]]$. Here, $\rho^L$ and $\rho^R$ are also extended by $\mathbb{C}[[t]]$-bilinearity.
%\begin{eqnarray*}
%&&\rho^L(\sum_{i=0}^{+\infty}x_it^{i})(\sum_{j=0}^{+\infty}v_jt^{j})=\sum_{k=0}^{+\infty}\sum_{i+j=k}\rho^L(x_i)v_jt^{k},\\
%&&\rho^R(\sum_{i=0}^{+\infty}x_it^{i})(\sum_{j=0}^{+\infty}v_jt^{j})=\sum_{k=0}^{+\infty}\sum_{i+j=k}\rho^R(x_i)v_jt^{k}, ~~~~\forall x_i\in g, v_j\in V.
%\end{eqnarray*}
 Similarly, the $2$-cocycle $H$ can be extended to a $2$-cocycle (which we denote by the same notation $H$) on the Leibniz algebra $\mathfrak{g}[[t]]$ with coefficients in $V[[t]]$.

Let $K : V \rightarrow \mathfrak{g}$ be a $H$-twisted relative Rota-Baxter operator on the Leibniz algebra $(\mathfrak{g}, [\c, \c])$ with respect to the representation $(V, \rho^L, \rho^{R})$ and $2$-cocycle $H$.  We consider the  power series
\begin{eqnarray*}
K_t=\sum_{i=0}^{+\infty}K_it^{i},~ \text{ for } K_i\in Hom (V, \mathfrak{g}) ~ \text{ with } K_0 = K.
\end{eqnarray*}
Extend $K_t$ to a linear map from $V [[t]]$ to
$\mathfrak{g}[[t]]$ by $\mathbb{C}[[t]]$-linearity which we still denote by $K_t$.

\begin{definition}
A formal deformation of $K$ is given by a formal sum $K_t=\sum_{i=0}^{+\infty}K_it^{i}$ with $K_0 = K$ satisfying
\begin{align}\label{deform-eqn}
&&[K_tu, K_tv]=K_t \big( \rho^{L}(K_tu)v+\rho^{R}(K_tv)u+H(K_t(u), K_t(v)) \big), ~ \text{ for } u, v \in V.
\end{align}
\end{definition}
It follows that $K_t$ is a $H$-twisted relative Rota-Baxter operator on the Leibniz algebra $\mathfrak{g}[[t]]$ with respect to the representation $V[[t]]$ and $2$-cocycle $H$.

\begin{remark}
If $K_t=\sum_{i=0}^{+\infty}K_it^{i}$ is a formal deformation of a $H$-twisted relative Rota-Baxter operator $K$ on a Leibniz algebra $(\mathfrak{g}, [\c, \c])$  with respect to a representation $(V, \rho^L, \rho^{R})$ and $2$-cocycle $H$, then $[\c, \c]_{K_t}$ defined by
\begin{eqnarray*}
[u, v]_{K_t}:=\sum_{i=0}^{+\infty} \big( \rho^{L}(K_iu)v+\rho^{R}(K_iv)u+ \sum_{j+k = i} H(K_ju, K_k v) \big)t^{i},~ \text{ for } u, v\in V,
\end{eqnarray*}
is a formal deformation of the associated Leibniz algebra $(V,   [\c, \c]_{K})$.
\end{remark}

By expanding the identity (\ref{deform-eqn}) and comparing coefficients of various powers of $t$, we obtain for $n \geq 0$,
\begin{align*}
 \sum_{i+j=n}[K_i u, K_j v] = \sum_{i+j=n}K_i(\rho^{L}(K_ju)v+\rho^{R}(K_jv)u)+\sum_{i+j+k=n}K_iH(K_j(u), K_k(v)),
\end{align*}
for $u, v \in V$. It holds for $n = 0$ as $K$ is a $H$-twisted relative Rota-Baxter operator. For $n = 1$, we obtain
\begin{align*}
 [Ku,  K_1v]+[K_1u,  Kv] =~& K_1(\rho^{L}(Ku)v+\rho^{R}(Kv)u+H(Ku, Kv))\\
 &+K(\rho^{L}(K_1u)v+\rho^{R}(K_1v)u+H(K_1(u), Kv)+H(K(u), K_1v)).
\end{align*}
This condition is equivalent to $(\partial_K (K_1))(u, v) = 0$, for $u, v \in V$.

\begin{proposition}\label{linear-term-cocycle}
Let $K_t=\sum_{i=0}^{+\infty}K_it^{i}$ be a formal deformation of a $H$-twisted relative Rota-Baxter operator $K$. Then $K_1$ is a $1$-cocycle in the cohomology of the $H$-twisted relative Rota-Baxter operator
$K$, that is, $\partial_{K} (K_1)=0$.
\end{proposition}

\begin{definition}
The $1$-cocycle $K_1$ is called the infinitesimal of the formal deformation $K_t = \sum_{i=0}^{+\infty}K_it^{i}$.
\end{definition}

Next, we define an equivalence between two formal deformations of a $H$-twisted relative Rota-Baxter operator.

\begin{definition}
Two formal deformations $K_t = =\sum_{i=0}^{+\infty}K_it^{i}$ and $K'_t = =\sum_{i=0}^{+\infty}K'_it^{i}$ of a $H$-twisted relative Rota-Baxter operator $K$ are said to be equivalent if there exists an element $x \in \mathfrak{g}$, linear maps
$\phi_i \in gl(\mathfrak{g})$ and $\psi_i \in gl(V)$ for $i \geq 2$ such that the pair
\begin{align*}
\big(\phi_t =Id_\mathfrak{g} + t L_x + \sum_{i=2}^{+\infty}\phi_i t^{i},~ \psi_t=Id_V + t(\rho^{L}(x)+H(x, K-)) + \sum_{i=2}^{+\infty}\psi_i t^{i} \big)
\end{align*}
is a morphism of $H$-twisted relative Rota-Baxter operators from $K_t$ to $K'_t$.
\end{definition}
By equating coefficients of $t$ from
both sides of the identity $\phi_t \circ K_t = K'_t\circ \psi_t$, we obtain
\begin{eqnarray*}
&&K_1(u)-K'_1(u)=K(\rho^{L}(x)u+H(x, Ku))-[x, Ku]=d_{K}(x)(u), \text{ for } u \in V.
\end{eqnarray*}

As a summary, we get the following.

\begin{theorem}
The infinitesimal of a formal deformation of a $H$-twisted  relative Rota-Baxter operator $K$ is a $1$-cocycle in the cohomology of $K$, and
the corresponding cohomology class depends only on the equivalence class of the deformation of $K$.
\end{theorem}

\begin{definition}
A $H$-twisted relative Rota-Baxter operator $K$ is said to be rigid if any formal deformation of $K$ is
equivalent to the undeformed deformation $K'_t = K$ .
\end{definition}

In the next theorem, we give a sufficient condition for the rigidity of a $H$-twisted relative Rota-Baxter operator in terms of Nijenhuis elements.

\begin{theorem}
Let $K$ be a $H$-twisted relative Rota-Baxter operator. If $Z^1_K (V, \mathfrak{g}) = \partial_K (\mathrm{Nij}(K ))$ then $K$ is rigid.
\end{theorem}

{\bf Proof.} Let $K_t=\sum_{i=0}^{+\infty}K_it^{i}$ be any formal deformation of $K$. Then it follows from Proposition \ref{linear-term-cocycle} that the linear term $K_1$ is a $1$-cocycle in the
cohomology of $K$, i.e., $K_1 \in Z^1_K (V, \mathfrak{g})$. Thus, by the hypothesis, there is a Nijenhuis element $x \in \mathrm{Nij}(K)$
such that $K_1 = - \partial_K(x)$. We take
\begin{eqnarray*}
&& \phi_t =Id_\mathfrak{g} + tL_x ~~ \text{ and } ~~\psi_t = Id_V + t(\rho^{L}(x)+H(x, K-)),
\end{eqnarray*}
and define $K'_t= \phi_t \circ K_t \circ \psi_t^{-1}$. Then $K'_t$ is a formal deformation equivalent to $K_t$. For $u \in V$, we observe that
%\begin{eqnarray*}
% K'_t(u)&=&(Id_g-L_xt+L_x^{2}t^2+\c \c \c +(-1)^{i}L_x^{i}t^i)(K_t(u+\rho^{L}(x)u+H(x,Ku)))\\
%&=& K(u)+(K_1(u)+ K( \rho^L(x)(u)+H(x, Ku)) - [x, Ku]_g)t+K'_2(u)t^2+\c \c \c\\
%&=& K(u)+K'_2(u)t^2+\c \c \c.
%\end{eqnarray*}
\begin{align*}
 K'_t(u)=~&(Id_\mathfrak{g} + tL_x)(K_t(u - t \rho^{L}(x)u- tH(x,Ku) + \text{ power of } t^{ \geq 2}))\\
=~& K(u) + t (K_1 u - K \rho^L (x) u - KH (x, Ku) + [x, Ku])+ \text{ power of } t^{\geq 2}. \\
=~& K(u) + t^2 K'_2 (u) + \c \c \c \qquad (\text{as } K_1 = - \partial_K (x)).
\end{align*}
Hence the coefficient of $t$ in the expression of $K_t'$ is trivial. Applying the same process repeatedly, we get that $K_t$ is equivalent to $K$. Therefore,
$K$ is rigid. \hfill $\square$

\section{NS-Leibniz algebras}
\def\theequation{\arabic{section}. \arabic{equation}}
\setcounter{equation} {0}
In this section, we introduce NS-Leibniz algebras as the underlying structure of twisted Rota-Baxter operators. Here we study some properties of NS-Leibniz algebras and give some examples. Further study on NS-Leibniz algebras is postponed to a forthcoming paper.

\begin{definition}\label{defn-ns}
An NS-Leibniz algebra is a quadruple $(A,  \triangleright, \triangleleft,  \diamond )$ consisting of a  vector space $A$ together with three bilinear operations $\triangleright, \triangleleft, \diamond : A \o A \rightarrow  A$   satisfying for all $x, y,z \in A$,
\begin{eqnarray*}
&&(A1)~~ x \triangleright (y \ast z) = (x \triangleright y) \triangleright z + y \triangleleft (x \triangleright z), \\
&&(A2)~~ x\triangleleft(y\triangleright z)=(x\triangleleft y)\triangleright z+y\triangleright (x\ast z),\\
&&(A3)~~ x\triangleleft (y\triangleleft z) = (x \ast y)\triangleleft z  + y\triangleleft (x\triangleleft z),\\
&&(A4)~~ x \triangleleft (y  \diamond  z) + x \diamond  (y\ast z)  = (x   \diamond  y) \triangleright z + (x \ast y) \diamond  z  + y \triangleleft (x \diamond  z) +  y \diamond  (x \ast z) ,
\end{eqnarray*}
where $x \ast y = x \triangleright y + x \triangleleft y + x \diamond y$.
\end{definition}

NS-Leibniz algebras are more general than Leibniz-dendriform algebras introduced in \cite{ST19}. More precisely, an NS-Leibniz algebra $(A,  \triangleright, \triangleleft,  \diamond )$  in which the bilinear operation $\diamond$ is trivial is a Leibniz-dendriform algebra.

In the following, we show that NS-Leibniz algebras split Leibniz algebras.

\begin{proposition}
Let $(A, \triangleright, \triangleleft,  \diamond )$ be an NS-Leibniz algebra. Then the vector space $A$ with the bilinear operation
\begin{align*}
[\c, \c]_\ast : A \o A \rightarrow A, ~[x, y]_\ast := x \ast y
\end{align*}
is a Leibniz algebra.
\end{proposition}
{\bf Proof.} By summing up the left hand sides of the identities (A1)-(A4), we simply get $[x , [ y , z]_\ast ]_\ast$. On the other hand, by summing up the right hand sides of the identities (A1)-(A4), we get $[[x , y]_\ast , z]_\ast + [y , [x , z]_\ast ]_\ast.$ Hence the result follows.
\hfill $\square$

\medskip

The Leibniz algebra $(A, [\c, \c]_\ast)$ of the above proposition is called the subadjacent Leibniz algebra of $(A,  \triangleright, \triangleleft, \diamond)$ and
$(A, \triangleright, \triangleleft, \diamond)$ is called a compatible NS-Leibniz algebra structure on $(A, [\c, \c]_{\ast}).$

\begin{proposition}
Let $(\mathfrak{g}, [\c, \c])$ be a Leibniz algebra and $N: \mathfrak{g} \rightarrow \mathfrak{g}$ be a Nijenhuis operator on it. Then the bilinear operations
\begin{eqnarray*}
x \triangleright y = [x, Ny],~~ x \triangleleft y = [Nx, y] ~~ \text{ and } ~~ x \diamond y = -N[x, y],~\text{ for } x, y\in \mathfrak{g}
\end{eqnarray*}
defines an NS-Leibniz algebra structure on $\mathfrak{g}$.
\end{proposition}

{\bf Proof.} For any $x, y, z\in \mathfrak{g}$, we have
\begin{align*}
    x \triangleright (y \ast z) = [x, N (y \ast z)]
    =~& [x, [Ny, Nz]] \\
    =~& [[x, Ny], Nz] + [Ny, [x, Nz]] \\
    =~& (x \triangleright y ) \triangleright z + y \triangleleft (x \triangleright z).
\end{align*}
Hence the identity (A1) of Definition \ref{defn-ns} holds. Similarly,
\begin{align*}
    x \triangleleft (y \triangleright z) = [Nx, [y, Nz]]
    =~& [[Nx, y], Nz] + [y, [Nx, Nz]] \\
    =~& (x \triangleleft y ) \triangleright z + y \triangleright (x \ast z),
\end{align*}
and
\begin{align*}
    x \triangleleft ( y \triangleleft z ) = [Nx, [Ny,z] ]
    =~& [[Nx, Ny], z] + [Ny, [Nx, z]] \\
    =~& (x \ast y) \triangleleft z + y \triangleleft ( x \triangleleft z).
\end{align*}
Therefore, the identities (A2) and (A3) also hold. To prove the identity (A4), we first recall from \cite{grab} that the given Leibniz bracket $[\c, \c]$ and the deformed Leibniz bracket $[ \c, \c]_N$ given in (\ref{deformed-leib}) are compatible in sense that their sum also defines a Leibniz bracket on $\mathfrak{g}$. This is equivalent to the fact that
\begin{align}\label{comp-equiv}
    [x,[y,z]]_N   +  [x,[y,z]_N]  =  [[x,y],z]_N + [[x,y]_N,z] + [y, [x,z]]_N + [y, [x, z]_N],
\end{align}
for $x, y, z \in \mathfrak{g}$. The identity (A4) of Definition \ref{defn-ns} simply follows from (\ref{comp-equiv}). Hence $(\mathfrak{g}, \triangleright, \triangleleft, \diamond)$ is an NS-Leibniz algebra. \hfill $\square$

\medskip

Let $(A,  \triangleright, \triangleleft, \diamond)$ be an NS-Leibniz algebra. Define two linear maps $L_{\triangleleft}: A \rightarrow gl(A)$,~
$R_{\triangleright}: A \rightarrow gl(A)$ and a bilinear map $H : A \otimes A \rightarrow A$ by
\begin{align*}
L_{\triangleleft}(x)y=x\triangleleft y,~~~~~~~ R_{\triangleright}(x)y=y\triangleright x,~~~~~~H(x, y)=x \diamond y,~ \text{ for } x, y\in A.
\end{align*}
With these notations, we have the following.

\begin{proposition}\label{ns-compatible}
Let $(A, \triangleright, \triangleleft, \diamond)$ be an NS-Leibniz algebra. Then $(A,  L_{\triangleleft}, R_{\triangleright})$ is a representation of the subadjacent Leibniz algebra $(A, [\c, \c]_{\ast})$, and $H$ defined above is a $2$-cocycle. Moreover, the identity map $Id: A \rightarrow A$ is a $H$-twisted relative Rota-Baxter operator on the Leibniz algebra $(A, [\c, \c]_{\ast})$ with respect to the representation $(A,  L_{\triangleleft}, R_{\triangleright})$.
\end{proposition}

{\bf Proof.} For any $x, y, z\in A$, we have
\begin{eqnarray*}
 L_{\triangleleft}([x, y]_{\ast})z = [x, y]_{\ast} \triangleleft z
&\stackrel{(A3)}{=}& x\triangleleft (y\triangleleft z)-y\triangleleft(x\triangleleft z)\\
&=& \big( L_{\triangleleft}(x)\circ L_{\triangleleft}(y)-L_{\triangleleft}(y)\circ L_{\triangleleft}(x) \big) z.
\end{eqnarray*}
Similarly,
\begin{eqnarray*}
R_{\triangleright}([x,y]_\ast) z = z\triangleright [x, y]_{\ast}
&\stackrel{(A2)}{=}& x\triangleleft (z\triangleright y)-(x\triangleleft z)\triangleright y\\
&=&L_{\triangleleft}(x)\circ R_{\triangleright}(y))z-R_{\triangleright}(y)\circ L_{\triangleleft}(x)z,
\end{eqnarray*}
and
\begin{eqnarray*}
R_\triangleright ([x,y]_*)z
= z \triangleright [x,y]_\ast
&\stackrel{(A1)}{=}& (z \triangleright x) \triangleright y + x \triangleleft (z \triangleright y) \\
&=& \big( R_{\triangleright}(y) \circ R_{\triangleright}(x) + L_{\triangleleft}(x)R_{\triangleright}(y) \big) z.
\end{eqnarray*}
Therefore, $(A,  L_{\triangleleft}, R_{\triangleright})$ is a representation of the subadjacent Leibniz algebra $(A, [\c, \c]_{\ast})$. Moreover, the condition (A4) is equivalent that $H$ is a $2$-cocycle in the Loday-Pirashvili cochain complex of the Leibniz algebra $(A, [\c, \c]_\ast)$ with coefficients in the representation $(A, L_\triangleleft, R_\triangleright)$. Finally, we have
\begin{eqnarray*}
Id(L_{\triangleleft}(Id ~x)y+R_{\triangleright}(Id ~y)x+H(Id~x, Id ~y))=x\triangleleft y+x\triangleright y+x\circ y=[Id ~x, Id ~y]_{\ast},
\end{eqnarray*}
which shows that $Id : A \rightarrow A$ is a $H$-twisted relative Rota-Baxter operator on the Leibniz algebra
$(A, [\c, \c]_{\ast})$ with respect to the representation $(A,  L_{\triangleleft}, R_{\triangleright})$.  \hfill $\square$

\begin{proposition}\label{twisted-rota-ns}
Let $(\mathfrak{g}, [\c, \c])$ be a Leibniz algebra, $(V, \rho^L, \rho^{R})$ be a representation and $H \in C^2 (\mathfrak{g}, V)$ be a $2$-cocycle. Let $K : V \rightarrow \mathfrak{g}$ be a $H$-twisted relative Rota-Baxter operator. Then there is an NS-Leibniz algebra structure on $V$ with bilinear operations given by
\begin{eqnarray*}
u \triangleright v := \rho^{R}(Kv)u,~~~~ u \triangleleft v := \rho^{L}(Ku)v,~~~ u \diamond v := H(Ku, Kv), ~~\text{ for } u, v\in V.
\end{eqnarray*}
\end{proposition}

{\bf Proof.} For any $u, v, w\in V$, we have
\begin{align*}
u \triangleright (v \ast w) = \rho^R ( K (v \ast w)) u
=~& \rho^R ([Kv, Kw]) u \\
=~& \rho^L (Kv) \rho^R (Kw) u + \rho^R (Kw) \rho^R (Kv) u \\
=~& v \triangleleft (u \triangleright w ) + ( u \triangleright v) \triangleright w.
\end{align*}
Similarly,
\begin{align*}
    u \triangleleft (v \triangleright w) = \rho^L (Ku) \rho^R (Kw) v
    =~& \rho^R ([Ku, Kz])v + \rho^R (Kw) \rho^L (Ku) v \\
    =~& v \triangleright (u \ast w) + ( u \triangleleft v) \triangleright w,
\end{align*}
and
\begin{align*}
    u \triangleleft ( v \triangleleft w ) = \rho^L (Ku) \rho^L (Kv) (w)
    =~& \rho^L ([Ku, Kv]) w + \rho^L (Kv) \rho^L (Ku) w \\
    =~& (u \ast v) \triangleleft w +  v \triangleleft ( u \triangleleft w).
\end{align*}
Hence (A1), (A2) and (A3) of Definition \ref{defn-ns} holds. Since $H$ is a $2$-cocycle, we have $(\partial H)(Ku, Kv, Kz) = 0$, i.e.,
\begin{align*}
    \rho^{L}(Ku)H(Kv, Kw) &- \rho^{L}(Kv)H(Ku, Kw) - \rho^{R}(Kw)H(Ku, Kv) \\
    -& H([Ku, Kv], Kw) - H(Kv, [Ku, Kw]) + H(Ku, [Kv, Kw])=0.
\end{align*}
This is equivalent to the condition (A4) of Definition \ref{defn-ns}. Hence the proof.
\hfill $\square$

\begin{remark}
The subadjacent Leibniz algebra of the NS-Leibniz algebra constructed in Proposition \ref{twisted-rota-ns} is given by
\begin{align*}
    [u, v]_* = \rho^L (Ku) v + \rho^R (Kv) u + H (Ku, Kv), ~ \text{ for } u, v \in V.
\end{align*}
This Leibniz algebra on $V$ coincides with the one given in Proposition \ref{induced-leib}.
\end{remark}

\medskip

In the following, we give a necessary and sufficient condition for the existence of a compatible NS-Leibniz algebra structure on a Leibniz algebra.

\begin{proposition}
Let $(\mathfrak{g}, [\c, \c])$ be a Leibniz algebra. Then there is a compatible NS-Leibniz algebra structure on $\mathfrak{g}$ if and only if there exists an invertible $H$-twisted relative Rota-Baxter operator $K: V\rightarrow \mathfrak{g}$ on
$\mathfrak{g}$ with respect to a representation $(V, \rho^L, \rho^{R})$ and a $2$-cocycle $H$. Furthermore, the compatible NS-Leibniz algebra structure on $\mathfrak{g}$ is given by
\begin{align*}
x\triangleright y := K(\rho^{R}(y)K^{-1}x),~~~ x\triangleleft y:=K(\rho^{L}(x)K^{-1}y),~~~ x \diamond y=KH(x,y),~ \text{ for } x, y\in \mathfrak{g}.
\end{align*}
\end{proposition}

{\bf Proof.}  Let $K: V\rightarrow \mathfrak{g}$ be an invertible $H$-twisted relative Rota-Baxter operator on $\mathfrak{g}$ with respect to a representation $(V,  \rho^L, \rho^{R})$ and a $2$-cocycle $H$. By Proposition \ref{twisted-rota-ns}, there is an NS-Leibniz algebra structure on $V$ given by
\begin{eqnarray*}
u ~\Bar{\triangleright}~ v:=\rho^{R}(Kv)u,~~~~~u ~\Bar{\triangleleft}~ v:=\rho^{L}(Ku)v,~~~~~u ~\Bar{\diamond}~ v:=H(Ku, Kv),~ \text{ for } u, v\in V.
\end{eqnarray*}
Since $K$ is an invertible map, the bilinear operations
\begin{eqnarray*}
&&x \triangleright y := K(K^{-1} x ~\Bar{\triangleright}~ K^{-1} y) = K(\rho^{R}(y)K^{-1}x), \\
&&x \triangleleft y := K(K^{-1} x ~\Bar{\triangleleft}~ K^{-1}y) = K(\rho^{L}(x)K^{-1}y),\\
&& x \diamond y := K(K^{-1} x ~\Bar{\diamond}~ K^{-1}y) =KH(x,y), ~\text{ for } x, y \in \mathfrak{g}
\end{eqnarray*}
defines an NS-Leibniz algebra on $\mathfrak{g}$. Moreover,  we have
\begin{eqnarray*}
 && x \triangleright y+x \triangleleft y+x\diamond y\\
 &=& K(\rho^{R}(y)K^{-1}x)+K(\rho^{L}(x)K^{-1}y)+KH(x,y) \\
 &=& K(\rho^{R}(K\circ K^{-1}y)K^{-1}x)+K(\rho^{L}(K\circ K^{-1}x)K^{-1}y)+KH(K\circ K^{-1}x,K\circ K^{-1}y)\\
 &=&[K\circ K^{-1}x, K\circ K^{-1}y] = [x,y].
\end{eqnarray*}

Conversely, let $(\mathfrak{g}, \triangleright, \triangleleft, \diamond)$ be a compatible NS-Leibniz algebra structure on $\mathfrak{g}$. By Proposition \ref{ns-compatible},  $(\mathfrak{g}, L_{\triangleleft}, R_{\triangleright})$ is a representation of the Leibniz algebra $(\mathfrak{g}, [\c, \c])$, and the identity map $Id: \mathfrak{g} \rightarrow \mathfrak{g}$ is a $H$-twisted relative Rota-Baxter operator on the Leibniz algebra $(\mathfrak{g}, [\c, \c])$ with
respect to the representation $(\mathfrak{g},   L_{\triangleleft}, R_{\triangleright})$. Hence the proof. \hfill $\square$

\medskip

\medskip

\begin{center}
 {\bf ACKNOWLEDGEMENT}
 \end{center}

The work of A. Das is supported by the fellowship of Indian Institute of Technology (IIT) Kanpur. The work of S. Guo is supported by the NSF of China (No. 11761017) and Guizhou Provincial  Science and Technology  Foundation (No. [2020]1Y005).

\renewcommand{\refname}{REFERENCES}

\end{document}